\documentclass[english,oneside]{amsart}
\usepackage[T1]{fontenc}
\usepackage[latin1]{inputenc}
\usepackage{amssymb}

\makeatletter
 \theoremstyle{plain}
\newtheorem{thm}{Theorem}[section]
  \theoremstyle{definition}
  \newtheorem{defn}[thm]{Definition}
  \theoremstyle{plain}
  \newtheorem{cor}[thm]{Corollary}
  \theoremstyle{remark}
  \newtheorem{rem}[thm]{Remark}
  \theoremstyle{plain}
  \newtheorem{lem}[thm]{Lemma}

\usepackage{babel}
\makeatother

\begin{document}

\title{Scaled Asymptotics For Some $q$-Series As $q$ Approaching Unit }

\author{Ruiming Zhang}


\email{ruimingzhang@yahoo.com}

\begin{abstract}
In this work we investigate Plancherel-Rotach type asymptotics for
some $q$-series as $q\to1$. These $q$-series generalize Ramanujan
function $A_{q}(z)$; Jackson's $q$-Bessel function $J_{\nu}^{(2)}$(z;q),
Ismail-Masson orthogonal polynomials($q^{-1}$-Hermite polynomials)
$h_{n}(x|q)$, Stieltjes-Wigert orthogonal polynomials $S_{n}(x;q)$,
$q$-Laguerre orthogonal polynomials $L_{n}^{(\alpha)}(x;q)$ and
confluent basic hypergeometric series. 
\end{abstract}

\curraddr{School of Mathematics\\
Guangxi Normal University\\
Guilin City, Guangxi 541004\\
P. R. China.}


\maketitle

\section{Introduction\label{sec:Introduction}}

In \cite{Zhang1} we derived certain Plancherel-Rotach type asymptotics
for some $q$-series. These $q$-series generalize Ramanujan's entire
function $A_{q}(z)$, Jackson's $q$-Bessel function $J_{\nu}^{(2)}$(z;q),
Ismail-Masson orthogonal polynomials ($q^{-1}$-Hermite polynomials)
$h_{n}(x|q)$, Stieltjes-Wigert orthogonal polynomials $S_{n}(x;q)$,
$q$-Laguerre orthogonal polynomials $L_{n}^{(\alpha)}(x;q)$ and
confluent basic hypergeometric series. 

In this work we shall employ the method used in \cite{Zhang2} to
study the scaled asymptotics of these $q$-series as $q\to1$. In
section \S\ref{sec:Preliminaries} we list some common notations
from $q$-series and special functions. We present our results in
section \S\ref{sec:Main-Results} and prove them in section \S\ref{sec:Proofs}.
Throughout this work we always assume that $0<q<1$ unless otherwise
stated.

\section{Preliminaries\label{sec:Preliminaries}}

For a complex number $z$, we define \cite{Andrews,Gasper,Ismail2,Koekoek}\begin{equation}
(z;q)_{\infty}:=\prod_{k=0}^{\infty}(1-zq^{k}),\label{eq:1}\end{equation}
and the $q$-Gamma function is defined as\begin{equation}
\Gamma_{q}(z):=\frac{(q;q)_{\infty}}{(q^{z};q)_{\infty}}(1-q)^{1-z}\quad z\in\mathbb{C}.\label{eq:2}\end{equation}
The $q$-shifted factorials of $a,a_{1},\dotsc a_{m}$ are given by

\begin{equation}
(a;q)_{n}:=\frac{(a;q)_{\infty}}{(aq^{n};q)_{\infty}},\quad(a_{1},\dotsc,a_{m};q)_{n}:=\prod_{k=1}^{m}(a_{k};q)_{n}\label{eq:3}\end{equation}
for all integers $n\in\mathbb{Z}$ and $m\in\mathbb{N}$. Given two
nonnegative integers $s,r$ and two sets of complex numbers $a_{1},\dots,a_{r}$
and $b_{1},\dotsc,b_{s}$, a basic hypergeometric series ${}_{s}\phi_{r}$
is formally defined as\begin{equation}
{}_{s}\phi_{r}\left(\begin{array}{c}
a_{1},\dotsc,a_{r}\\
b_{1},\dotsc,b_{s}\end{array}\vert q,z\right):=\sum_{k=0}^{\infty}\frac{(a_{1},\dotsc,a_{r};q)_{k}(-zq^{-\ell})^{k}q^{\ell k^{2}}}{(q,b_{1},\dotsc,b_{s};q)_{k}},\label{eq:4}\end{equation}
where \begin{equation}
\ell:=\frac{s+1-r}{2},\label{eq:5}\end{equation}
and it is a confluent basic hypergeometric series if $\ell>0$. 

Given nonnegative integers $r,s,t$ and a positive number $\ell$,
we define \cite{Zhang1} \begin{align}
 & g(a_{1},\dotsc,a_{r};b_{1},\dotsc,b_{s};q;\ell;z)\label{eq:6}\\
 & :=\sum_{k=0}^{\infty}\frac{(q^{k+1},b_{1}q^{k},\dotsc,b_{s}q^{k};q)_{\infty}q^{\ell k^{2}}(-z)^{k}}{(a_{1}q^{k},\dotsc,a_{r}q^{k};q)_{\infty}},\nonumber \\
 & h(a_{1},\dotsc,a_{r};b_{1},\dotsc,b_{s};c_{1},\dots,c_{t};q;\ell;z)\label{eq:7}\\
 & :=\sum_{k=0}^{n}\frac{(q^{k+1},b_{1}q^{k},\dotsc,b_{s}q^{k};q)_{\infty}q^{\ell k^{2}}(-z)^{k}}{(a_{1}q^{k},\dotsc,a_{r}q^{k};q)_{\infty}}\frac{(q,c_{1},\dotsc,c_{t};q)_{n}}{(q,c_{1},\dotsc,c_{t};q)_{n-k}},\nonumber \end{align}
 where \begin{equation}
0\le a_{1},\dotsc,a_{r},b_{1},\dotsc,b_{s},c_{1},\dotsc,c_{t}<1.\label{eq:8}\end{equation}
Jackson's $q$-Bessel function $J_{\nu}^{(2)}(z;q)$ is defined as
\cite{Ismail2,Gasper,Koekoek} \begin{equation}
J_{\nu}^{(2)}(z;q):=\frac{(q^{\nu+1};q)_{\infty}}{(q;q)_{\infty}}\sum_{k=0}^{\infty}\frac{q^{k^{2}+k\nu}(-1)^{k}}{(q,q^{\nu+1};q)_{k}}\left(\frac{z}{2}\right)^{2k+\nu},\quad\nu>-1.\label{eq:9}\end{equation}
The Ismail-Masson polynomials $\left\{ h_{n}(x|q)\right\} _{n=0}^{\infty}$
are defined as \cite{Ismail2} \begin{equation}
h_{n}(\sinh\xi|q):=\sum_{k=0}^{n}\frac{(q;q)_{n}q^{k(k-n)}(-1)^{k}e^{(n-2k)\xi}}{(q;q)_{k}(q;q)_{n-k}}.\label{eq:10}\end{equation}
Stieltjes\emph{-}Wigert orthogonal polynomials $\left\{ S_{n}(x;q)\right\} _{n=0}^{\infty}$
are defined as \cite{Ismail2}\begin{equation}
S_{n}(x;q):=\sum_{k=0}^{n}\frac{q^{k^{2}}(-x)^{k}}{(q;q)_{k}(q;q)_{n-k}}.\label{eq:11}\end{equation}
The $q$-Laguerre orthogonal polynomials $\left\{ L_{n}^{(\alpha)}(x;q)\right\} _{n=0}^{\infty}$
are defined as \cite{Ismail2} \begin{equation}
L_{n}^{(\alpha)}(x;q):=\sum_{k=0}^{n}\frac{q^{k^{2}+\alpha k}(-x)^{k}(q^{\alpha+1};q)_{n}}{(q;q)_{k}(q,q^{\alpha+1};q)_{n-k}}\label{eq:12}\end{equation}
for $\alpha>-1$. Clearly, we have\begin{align}
A_{q}(z) & =\frac{g(-;-;q;1;z)}{(q;q)_{\infty}},\label{eq:13}\\
J_{\nu}^{(2)}(z;q) & =\frac{g(-;q^{\nu+1};q;1;z^{2}q^{\nu}/4)}{(q;q)_{\infty}^{2}(2/z)^{\nu}},\label{eq:14}\\
h_{n}(\sinh\xi\vert q) & =\frac{h(-;-;-;q;1;e^{-2\xi}q^{-n})}{e^{-n\xi}(q;q)_{\infty}},\label{eq:15}\\
S_{n}(x;q) & =\frac{h(-;-;-;q;1;x)}{(q;q)_{n}(q;q)_{\infty}},\label{eq:16}\\
L_{n}^{(\alpha)}(x;q) & =\frac{h(-;-;q^{\alpha+1};q;1;xq^{\alpha})}{(q;q)_{n}(q;q)_{\infty}},\label{eq:17}\\
{}_{s}\phi_{r}\left(\begin{array}{c}
a_{1},\dotsc,a_{r}\\
b_{1},\dotsc,b_{s}\end{array}\vert q,z\right) & =\frac{(q,b_{1},\dotsc,b_{s};q)_{\infty}g(a_{1},\dotsc,a_{r};b_{1},\dotsc b_{s};q;\ell;zq^{-\ell})}{(a_{1},\dotsc,a_{r};q)_{\infty}}.\label{eq:18}\end{align}
The Dedekind $\eta(\tau)$ is defined as \cite{Rademarcher}\begin{equation}
\eta(\tau):=e^{\pi i\tau/12}\prod_{k=1}^{\infty}(1-e^{2\pi ik\tau}),\label{eq:19}\end{equation}
 or\begin{equation}
\eta(\tau)=q^{1/12}(q^{2};q^{2})_{\infty},\quad q=e^{\pi i\tau},\quad\Im(\tau)>0.\label{eq:20}\end{equation}
 It has the transformation formula\begin{equation}
\eta\left(-\frac{1}{\tau}\right)=\sqrt{\frac{\tau}{i}}\eta(\tau).\label{eq:21}\end{equation}
The four Jacobi theta functions are defined as \cite{Rademarcher}
\begin{align}
\theta_{1}(v|\tau) & :=-i\sum_{k=-\infty}^{\infty}(-1)^{k}q^{(k+1/2)^{2}}e^{(2k+1)\pi iv},\label{eq:22}\\
\theta_{2}(v|\tau) & :=\sum_{k=-\infty}^{\infty}q^{(k+1/2)^{2}}e^{(2k+1)\pi iv},\label{eq:23}\\
\theta_{3}(v|\tau) & :=\sum_{k=-\infty}^{\infty}q^{k^{2}}e^{2k\pi iv},\label{eq:24}\\
\theta_{4}(v|\tau) & :=\sum_{k=-\infty}^{\infty}(-1)^{k}q^{k^{2}}e^{2k\pi iv},\label{eq:25}\end{align}
 where\begin{equation}
q=e^{\pi i\tau},\quad\Im(\tau)>0.\label{eq:26}\end{equation}
 For our convenience, we also use the following notations\begin{equation}
\theta_{\lambda}(z;q):=\theta_{\lambda}(v|\tau),\quad z=e^{2\pi iv},\quad q=e^{\pi i\tau}\label{eq:27}\end{equation}
 with \begin{equation}
\lambda=1,2,3,4.\label{eq:28}\end{equation}
 The Jacobi's triple product identities are\begin{align}
\theta_{1}(v|\tau) & =2q^{1/4}\sin\pi v(q^{2};q^{2})_{\infty}(q^{2}e^{2\pi iv};q^{2})_{\infty}(q^{2}e^{-2\pi iv};q^{2})_{\infty},\label{eq:29}\\
\theta_{2}(v|\tau) & =2q^{1/4}\cos\pi v(q^{2};q^{2})_{\infty}(-q^{2}e^{2\pi iv};q^{2})_{\infty}(-q^{2}e^{-2\pi iv};q^{2})_{\infty},\label{eq:30}\\
\theta_{3}(v|\tau) & =(q^{2};q^{2})_{\infty}(-qe^{2\pi iv};q^{2})_{\infty}(-qe^{-2\pi iv};q^{2})_{\infty},\label{eq:31}\\
\theta_{4}(v|\tau) & =(q^{2};q^{2})_{\infty}(qe^{2\pi iv};q^{2})_{\infty}(qe^{-2\pi iv};q^{2})_{\infty}.\label{eq:32}\end{align}
 The Jacobi $\theta$ functions satisfy transformations\begin{align}
\theta_{1}\left(\frac{v}{\tau}\mid-\frac{1}{\tau}\right) & =-i\sqrt{\frac{\tau}{i}}e^{\pi iv^{2}/\tau}\theta_{1}\left(v\mid\tau\right),\label{eq:33}\\
\theta_{2}\left(\frac{v}{\tau}\mid-\frac{1}{\tau}\right) & =\sqrt{\frac{\tau}{i}}e^{\pi iv^{2}/\tau}\theta_{4}\left(v\mid\tau\right),\label{eq:34}\\
\theta_{3}\left(\frac{v}{\tau}\mid-\frac{1}{\tau}\right) & =\sqrt{\frac{\tau}{i}}e^{\pi iv^{2}/\tau}\theta_{3}\left(v\mid\tau\right),\label{eq:35}\\
\theta_{4}\left(\frac{v}{\tau}\mid-\frac{1}{\tau}\right) & =\sqrt{\frac{\tau}{i}}e^{\pi iv^{2}/\tau}\theta_{2}\left(v\mid\tau\right).\label{eq:36}\end{align}
The Euler Gamma function $\Gamma(z)$ is given by \cite{Andrews,Gasper,Ismail2,Koekoek}
\begin{equation}
\frac{1}{\Gamma(z)}:=z\prod_{k=1}^{\infty}\left(1+\frac{z}{k}\right)\left(1+\frac{1}{k}\right)^{-z},\quad z\in\mathbb{C}.\label{eq:37}\end{equation}
For any real number $x$, we have \begin{equation}
x=\left\lfloor x\right\rfloor +\left\{ x\right\} ,\label{eq:38}\end{equation}
where the fractional part of $x$ is $\left\{ x\right\} \in[0,1)$
and $\left\lfloor x\right\rfloor \in\mathbb{Z}$ is the greatest integer
less or equal $x$. The arithmetic function \begin{equation}
\chi(n):=\begin{cases}
1 & 2\nmid n\\
0 & 2\mid n\end{cases},\label{eq:39}\end{equation}
 which is the principal character modulo $2$, satisfies the identities\begin{equation}
\chi(n)=2\left\{ \frac{n}{2}\right\} =n-2\left\lfloor \frac{n}{2}\right\rfloor =\left\lfloor \frac{n+1}{2}\right\rfloor -\left\lfloor \frac{n}{2}\right\rfloor .\label{eq:40}\end{equation}
 Thus,\begin{equation}
\left\lfloor \frac{n+1}{2}\right\rfloor =\frac{n+\chi(n)}{2},\label{eq:41}\end{equation}
 and\begin{equation}
\left\lfloor \frac{n}{2}\right\rfloor =\frac{n-\chi(n)}{2}.\label{eq:42}\end{equation}

\section{Main Results\label{sec:Main-Results}}

In order to state our results in full generality we also need the
following definition:

\begin{defn}
\label{def:admissible-scale}An admissible scale is a sequence $\left\{ \lambda_{n}\right\} _{n=1}^{\infty}$
of positive numbers such that\begin{equation}
\lim_{n\to\infty}\frac{\lambda_{n}}{\log n}=\infty,\quad\lim_{n\to\infty}\frac{n}{\lambda_{n}^{2}}=\infty.\label{eq:43}\end{equation}

\end{defn}
Clearly, \begin{equation}
\lambda_{n}=n^{\beta}\log^{\gamma}n,\quad\beta+\gamma>0,\quad0<\beta<\frac{1}{2},\quad\gamma\ge0,\label{eq:44}\end{equation}
and\begin{equation}
\lambda_{n}=\log^{\gamma}n,\quad\gamma>1\label{eq:45}\end{equation}
are admissible scales.

\subsection{$g$-function}

To simplify the type setting in the following theorem, we let 

\begin{equation}
g(z;q):=g(a_{1},\dotsc,a_{r};b_{1},\dotsc,b_{s};q;\ell;z).\label{eq:46}\end{equation}

\begin{thm}
\label{thm:g-asymptotics-1}Given an admissible scale $\lambda_{n}$,
assume that \begin{equation}
z=e^{2\pi v},\quad q:=e^{-\pi\lambda_{n}^{-1}},\quad\ell>0,\quad v\in\mathbb{R},\label{eq:47}\end{equation}
 and\begin{equation}
a_{j}:=q^{\alpha_{j}},\quad b_{k}:=q^{\beta_{k}},\quad\alpha_{j},\beta_{k}>0\label{eq:48}\end{equation}
for \begin{equation}
1\le j\le r,\quad1\le k\le s.\label{eq:49}\end{equation}
Then, \begin{align}
g(-q^{-4n\ell}z;q) & =\exp\left\{ \frac{\pi\lambda_{n}}{\ell}\left(v+\frac{2n\ell}{\lambda_{n}}\right)^{2}\right\} \label{eq:50}\\
 & \times\sqrt{\frac{\lambda_{n}}{\ell}}\left\{ 1+\mathcal{O}(e^{-\ell^{-1}\pi\lambda_{n}})\right\} ,\nonumber \end{align}
 and\begin{align}
g(q^{-4n\ell}z;q) & =\exp\left\{ \frac{\pi\lambda_{n}}{\ell}\left(v+\frac{2n\ell}{\lambda_{n}}\right)^{2}-\frac{\pi\lambda_{n}}{4\ell}\right\} \label{eq:51}\\
 & \times2\sqrt{\frac{\lambda_{n}}{\ell}}\left\{ \cos\frac{\pi\lambda_{n}v}{\ell}+\mathcal{O}\left(e^{-2\ell^{-1}\pi\lambda_{n}}\right)\right\} \nonumber \end{align}
 as $n\to\infty$, and the $\mathcal{O}$-term is uniform for $v$
in any compact subset of $\mathbb{R}$. 
\end{thm}
For the Ramanujan's entire function we have:

\begin{cor}
\label{cor:ramanujan}Given an admissible scale $\lambda_{n}$ , assume
that\begin{equation}
z=e^{2\pi v},\quad q=e^{-\pi\lambda_{n}^{-1}},\quad v\in\mathbb{R},\label{eq:52}\end{equation}
we have \begin{align}
A_{q}(-q^{-4n}z) & =\exp\left\{ \pi\lambda_{n}\left(v+\frac{2n}{\lambda_{n}}\right)^{2}+\frac{\pi\lambda_{n}}{6}-\frac{\pi}{24\lambda_{n}}\right\} \label{eq:54}\\
 & \times\frac{1}{\sqrt{2}}\left\{ 1+\mathcal{O}\left(e^{-\pi\lambda_{n}}\right)\right\} ,\nonumber \end{align}
 and\begin{align}
A_{q}(q^{-4n}z) & =\exp\left\{ \pi\lambda_{n}\left(v+\frac{2n}{\lambda_{n}}\right)^{2}-\frac{\pi\lambda_{n}}{12}-\frac{\pi}{24\lambda_{n}}\right\} \label{eq:55}\\
 & \times\sqrt{2}\left\{ \cos\pi\lambda_{n}v+\mathcal{O}\left(e^{-2\pi\lambda_{n}}\right)\right\} \nonumber \end{align}
as $n\to\infty$, and the $\mathcal{O}$-term is uniform for $v$
in any compact subset of $\mathbb{R}$. 
\end{cor}
For the Jackson's $q$-Bessel function we have:

\begin{cor}
\label{cor:jackson}For an admissible scale $\lambda_{n}$, assume
that\begin{equation}
z=e^{2\pi v},\quad q=e^{-\pi\lambda_{n}^{-1}},\quad v\in\mathbb{R},\quad\nu>-1,\label{eq:55}\end{equation}
 then, 

\begin{align}
J_{\nu}^{(2)}(2i\sqrt{zq^{-\nu}}q^{-2n};q) & =\frac{\exp\left(\frac{\pi\lambda_{n}}{3}-\frac{\pi}{12\lambda_{n}}+\frac{\nu^{2}\pi}{4\lambda_{n}}+\frac{\nu\pi i}{2}\right)}{2\sqrt{\lambda_{n}}}\label{eq:56}\\
 & \times\exp\left\{ \pi\lambda_{n}\left(v+\frac{4n+\nu}{2\lambda_{n}}\right)^{2}\right\} \left\{ 1+\mathcal{O}(e^{-\pi\lambda_{n}})\right\} ,\nonumber \end{align}

and\begin{align}
J_{\nu}^{(2)}(2\sqrt{zq^{-\nu}}q^{-2n};q) & =\frac{\exp\left(\frac{\pi\lambda_{n}}{12}-\frac{\pi}{12\lambda_{n}}+\frac{\nu^{2}\pi}{4\lambda_{n}}\right)}{\sqrt{\lambda_{n}}}\label{eq:57}\\
 & \times\exp\left\{ \pi\lambda_{n}\left(v+\frac{4n+\nu}{2\lambda_{n}}\right)^{2}\right\} \left\{ \cos\pi\lambda_{n}v+\mathcal{O}\left(e^{-2\pi\lambda_{n}}\right)\right\} \nonumber \end{align}
 as $n\to\infty$, and the $\mathcal{O}$-term is uniform for $v$
in any compact subset of $\mathbb{R}$. 
\end{cor}
For the confluent basic hypergeometric series we have:

\begin{cor}
\label{cor:confluent}Given an admissible scale $\lambda_{n}$, assume
that

\begin{equation}
z=e^{2\pi v},\quad q=e^{-\pi\lambda_{n}^{-1}},\quad v\in\mathbb{R},\label{eq:58}\end{equation}
and\begin{equation}
\alpha_{j},\beta_{k}>0,\quad1\le j\le r,\quad1\le k\le s.\label{eq:59}\end{equation}
Let \begin{equation}
\ell:=\frac{s+1-r}{2}>0,\quad\rho:=\sum_{j=1}^{r}\alpha_{j}-\sum_{j=1}^{s}\beta_{j}+\ell-1.\label{eq:60}\end{equation}
 Then we have\begin{align}
{}_{s}\phi_{r}\left(\begin{array}{c}
q^{\alpha_{1}},\dotsc,q^{\alpha_{r}}\\
q^{\beta_{1}},\dotsc,q^{\beta_{s}}\end{array}\vert q,-zq^{-\ell(4n-1)}\right) & =\frac{2^{\ell}\pi^{\rho+\ell}\prod_{j=1}^{r}\Gamma(\alpha_{j})}{\sqrt{\ell}\lambda_{n}^{\rho-1/2}\exp\left(\ell\pi\lambda_{n}/3\right)\prod_{j=1}^{s}\Gamma(\beta_{j})}\label{eq:61}\\
 & \times\left\{ \exp\frac{\pi\lambda_{n}}{\ell}\left(v+\frac{2n\ell}{\lambda_{n}}\right)^{2}\right\} \left\{ 1+\mathcal{O}\left(\lambda_{n}^{-1}\log^{2}\lambda_{n}\right)\right\} ,\nonumber \end{align}
 and\begin{align}
{}_{s}\phi_{r}\left(\begin{array}{c}
q^{\alpha_{1}},\dotsc,q^{\alpha_{r}}\\
q^{\beta_{1}},\dotsc,q^{\beta_{s}}\end{array}\vert q,zq^{-\ell(4n-1)}\right) & =\frac{2^{\ell+1}\pi^{\rho+\ell}\prod_{j=1}^{r}\Gamma(\alpha_{j})}{\sqrt{\ell}\lambda_{n}^{\rho-1/2}\exp\left(\frac{\ell\pi\lambda_{n}}{3}+\frac{\pi\lambda_{n}}{4\ell}\right)\prod_{j=1}^{s}\beta_{j}}\label{eq:62}\\
 & \times\left\{ \exp\frac{\pi\lambda_{n}}{\ell}\left(v+\frac{2n\ell}{\lambda_{n}}\right)^{2}\right\} \left\{ \cos\frac{\pi\lambda_{n}v}{\ell}+\mathcal{O}\left(\lambda_{n}^{-1}\log^{2}\lambda_{n}\right)\right\} \nonumber \end{align}
 as $n\to\infty$, and the $\mathcal{O}$-term is uniform for $v$
in any compact subset of $\mathbb{R}$. 
\end{cor}

\subsection{$h$-function}

For our convenience we let\begin{equation}
h_{n}(z;q):=h_{\ell}(a_{1},\dotsc,a_{r};b_{1},\dotsc,b_{s};c_{1},\dots,c_{t};q;z).\label{eq:63}\end{equation}
We have similar results for the $h$ function:

\begin{thm}
\label{thm:h-asymptotics-1}Given an admissible scale $\lambda_{n}$
, assume that \begin{equation}
z:=e^{2\pi v},\quad q:=e^{-\pi\lambda_{n}^{-1}},\quad\ell>0,\quad v\in\mathbb{R},\label{eq:64}\end{equation}
 and\begin{equation}
a_{j}:=q^{\alpha_{j}},\quad b_{k}:=q^{\beta_{k}},\quad\alpha_{j},\beta_{k}>0\label{eq:65}\end{equation}
for \begin{equation}
1\le j\le r,\quad1\le k\le s.\label{eq:66}\end{equation}
Then, \begin{align}
h(-zq^{-n\ell};q) & =\exp\left\{ \frac{\pi\lambda_{n}}{\ell}\left(v+\frac{\ell(n-\chi(n))}{2\lambda_{n}}\right)^{2}+\frac{\ell\pi(n-1)\chi(n)}{2\lambda_{n}}\right\} \label{eq:67}\\
 & \times\sqrt{\frac{\lambda_{n}}{\ell}}\left\{ 1+\mathcal{O}(e^{-\ell^{-1}\pi\lambda_{n}})\right\} ,\nonumber \end{align}
 and\begin{align}
h(zq^{-n\ell};q) & =\exp\left\{ \frac{\pi\lambda_{n}}{\ell}\left(v+\frac{\ell(n-\chi(n))}{2\lambda_{n}}\right)^{2}+\frac{\ell\pi(n-1)\chi(n)}{2\lambda_{n}}-\frac{\pi\lambda_{n}}{4\ell}\right\} \label{eq:68}\\
 & \times2\sqrt{\frac{\lambda_{n}}{\ell}}\left\{ \cos\frac{\pi\lambda_{n}}{\ell}\left(v+\frac{\ell(n-\chi(n))}{2\lambda_{n}}\right)+\mathcal{O}(e^{-2\ell^{-1}\pi\lambda_{n}})\right\} \nonumber \end{align}
as $n\to\infty$, and the $\mathcal{O}$-term is uniform for $v$
in any compact subset of $\mathbb{R}$. 
\end{thm}
For Ismail-Masson orthogonal polynomials we have:

\begin{cor}
\label{cor:ismail-masson} Given an admissible scale $\lambda_{n}$.
For any $v\in\mathbb{R}$ we have\begin{align}
h_{n}\left(\sinh\pi\left(v+\frac{i}{2}\right)\mid q\right) & =\frac{\exp\left\{ \frac{\pi n^{2}}{4\lambda_{n}}+\frac{\pi\lambda_{n}}{6}-\frac{\pi(1+12\chi(n))}{24\lambda_{n}}\right\} }{(-i)^{n}\sqrt{2}}\label{eq:69}\\
 & \times\left\{ \exp\left[\pi\lambda_{n}\left(v-\frac{\chi(n)}{2\lambda_{n}}\right)^{2}\right]\right\} \left\{ 1+\mathcal{O}(e^{-\pi\lambda_{n}})\right\} ,\nonumber \end{align}
 and\begin{align}
h_{n}(\sinh\pi v & \mid q)=(-1)^{n}\sqrt{2}\exp\left\{ \frac{n^{2}\pi}{4\lambda_{n}}-\frac{(1+12\chi(n))\pi}{24\lambda_{n}}-\frac{\pi\lambda_{n}}{12}\right\} \label{eq:70}\\
 & \times\left\{ \exp\left[\pi\lambda_{n}\left(v-\frac{\chi(n)}{2\lambda_{n}}\right)^{2}\right]\right\} \left\{ \cos\pi\lambda_{n}\left(v+\frac{n-\chi(n)}{2\lambda_{n}}\right)+\mathcal{O}(e^{-2\pi\lambda_{n}})\right\} \nonumber \end{align}

as $n\to\infty$, and the $\mathcal{O}$-term is uniform for $v$
in any compact subset of $\mathbb{R}$. 
\end{cor}
For Stieltjes-Wigert orthogonal polynomials we have:

\begin{cor}
\label{cor:stieltjes-wigert}Given an admissible scale $\lambda_{n}$,
assume that \begin{equation}
z:=e^{2\pi v},\quad q:=e^{-\pi\lambda_{n}^{-1}},\quad v\in\mathbb{R}.\label{eq:71}\end{equation}
Then we have\begin{align}
S_{n}(-zq^{-n};q) & =\frac{\exp\left\{ \frac{\pi\lambda_{n}}{3}+\frac{\pi(n-1)\chi(n)}{2\lambda_{n}}-\frac{\pi}{12\lambda_{n}}\right\} }{2\sqrt{\lambda_{n}}}\label{eq:72}\\
 & \times\left\{ \exp\pi\lambda_{n}\left(v+\frac{n-\chi(n)}{2\lambda_{n}}\right)^{2}\right\} \left\{ 1+\mathcal{O}\left(e^{-\pi\lambda_{n}}\right)\right\} ,\nonumber \end{align}
and\begin{align}
S_{n}(zq^{-n};q) & =\frac{\exp\left\{ \frac{\pi\lambda_{n}}{12}+\frac{\pi(n-1)\chi(n)}{2\lambda_{n}}-\frac{\pi}{12\lambda_{n}}\right\} }{\sqrt{\lambda_{n}}}\label{eq:73}\\
 & \times\left\{ \exp\pi\lambda_{n}\left(v+\frac{n-\chi(n)}{2\lambda_{n}}\right)^{2}\right\} \left\{ \cos\pi\lambda_{n}\left(v+\frac{n-\chi(n)}{2\lambda_{n}}\right)+\mathcal{O}\left(e^{-2\pi\lambda_{n}}\right)\right\} \nonumber \end{align}
as $n\to\infty$, and the $\mathcal{O}$-term is uniform for $v$
in any compact subset of $\mathbb{R}$. 
\end{cor}
For the $q$-Laguerre orthogonal polynomials we have:

\begin{cor}
\label{cor:laguerre}Given an admissible scale $\lambda_{n}$, assume
that\begin{equation}
z:=e^{-2\pi v},\quad q:=e^{-\pi\lambda_{n}^{-1}},\quad v\in\mathbb{R},\quad\alpha>-1.\label{eq:74}\end{equation}
 Then we have\begin{align}
L_{n}^{(\alpha)}(-zq^{-\alpha-n};q) & =\frac{\exp\left\{ \frac{\pi\lambda_{n}}{3}+\frac{\pi(n-1)\chi(n)}{2\lambda_{n}}-\frac{\pi}{12\lambda_{n}}\right\} }{2\sqrt{\lambda_{n}}}\label{eq:75}\\
 & \times\left\{ \exp\pi\lambda_{n}\left(v+\frac{n-\chi(n)}{2\lambda_{n}}\right)^{2}\right\} \left\{ 1+\mathcal{O}\left(e^{-\pi\lambda_{n}}\right)\right\} ,\nonumber \end{align}
 and\begin{align}
L_{n}^{(\alpha)}(zq^{-\alpha-n};q) & =\frac{\exp\left\{ \frac{\pi\lambda_{n}}{12}+\frac{\pi(n-1)\chi(n)}{2\lambda_{n}}-\frac{\pi}{12\lambda_{n}}\right\} }{\sqrt{\lambda_{n}}}\label{eq:76}\\
 & \times\left\{ \exp\pi\lambda_{n}\left(v+\frac{n-\chi(n)}{2\lambda_{n}}\right)^{2}\right\} \left\{ \cos\pi\lambda_{n}\left(v+\frac{n-\chi(n)}{2\lambda_{n}}\right)+\mathcal{O}\left(e^{-2\pi\lambda_{n}}\right)\right\} \nonumber \end{align}

as $n\to\infty$, and the $\mathcal{O}$-term is uniform for $v$
in any compact subset of $\mathbb{R}$. 
\end{cor}
\begin{rem}
Similar results hold for general $\tau$ and $\beta$ defined in \cite{Zhang1}
and their proofs are also similar to the proofs for the special cases
here. The formulas for the general $\tau$ and $\beta$ may be applicable
to the studies in phase transitions and critical phenomena in physics.
However, we feel that the formulas for the special cases are more
appealing and thus skip the general formuas. 
\end{rem}

\section{Proofs\label{sec:Proofs}}

The following lemma is from \cite{Zhang1}, we won't reproduce its
proof here.

\begin{lem}
\label{lem:1}Given a complex number $a$, assume that\begin{equation}
0<\frac{\left|a\right|q^{n}}{1-q}<\frac{1}{2}\label{eq:77}\end{equation}
 for some positive integer $n$. Then, \begin{equation}
\frac{(a;q)_{\infty}}{(a;q)_{n}}=(aq^{n};q)_{\infty}:=1+r_{1}(a;n)\label{eq:78}\end{equation}
 with\begin{equation}
\left|r_{1}(a;n)\right|\le\frac{2\left|a\right|q^{n}}{1-q}\label{eq:79}\end{equation}
 and \begin{equation}
\frac{(a;q)_{n}}{(a;q)_{\infty}}=\frac{1}{(aq^{n};q)_{\infty}}:=1+r_{2}(a;n)\label{eq:80}\end{equation}
 with \begin{equation}
\left|r_{2}(a;n)\right|\le\frac{2\left|a\right|q^{n}}{1-q}.\label{eq:81}\end{equation}

\end{lem}
We also need the following lemma:

\begin{lem}
\label{lem:2}Given a sequence of positive numbers $\left\{ \lambda_{n}\right\} _{n=1}^{\infty}$,
let \begin{equation}
q=e^{-\pi\lambda_{n}^{-1}},\quad\lim_{n\to\infty}\lambda_{n}=\infty,\label{eq:82}\end{equation}
 then,\begin{equation}
(q;q)_{\infty}=\sqrt{2\lambda_{n}}\exp\left(\frac{\pi}{24\lambda_{n}}-\frac{\pi\lambda_{n}}{6}\right)\left\{ 1+\mathcal{O}(e^{-4\pi\lambda_{n}})\right\} \label{eq:83}\end{equation}
 and\begin{equation}
\frac{1}{(q;q)_{\infty}}=\frac{\exp\left(\frac{\pi\lambda_{n}}{6}-\frac{\pi}{24\lambda_{n}}\right)}{\sqrt{2\lambda_{n}}}\left\{ 1+\mathcal{O}(e^{-4\pi\lambda_{n}})\right\} \label{eq:84}\end{equation}
 as $n\to\infty$. 
\end{lem}
\begin{proof}
From formulas \eqref{eq:19}, \eqref{eq:21} and \eqref{eq:21} we
get \begin{align}
 & (q;q)_{\infty}=\exp\left(\frac{\pi}{24\lambda_{n}}\right)\eta\left(\frac{i}{2\lambda_{n}}\right)\label{eq:85}\\
 & =\sqrt{2\lambda_{n}}\exp\left(\frac{\pi}{24\lambda_{n}}\right)\eta(2\lambda_{n}i)\nonumber \\
 & =\sqrt{2\lambda_{n}}\exp\left(\frac{\pi}{24\lambda_{n}}-\frac{\pi\lambda_{n}}{6}\right)\prod_{k=1}^{\infty}(1-e^{-4\pi k\lambda_{n}}).\nonumber \end{align}
For sufficiently large $n$ satisfying \begin{equation}
\exp(-4\pi\lambda_{n})<\frac{1}{3},\label{eq:86}\end{equation}
we have\begin{equation}
\prod_{k=1}^{\infty}(1-e^{-4\pi k\lambda_{n}})=1+\mathcal{O}(e^{-4\pi\lambda_{n}}),\label{eq:87}\end{equation}
and\begin{equation}
\frac{1}{\prod_{k=1}^{\infty}(1-e^{-4\pi k\lambda_{n}})}=1+\mathcal{O}(e^{-4\pi\lambda_{n}})\label{eq:88}\end{equation}
by Lemma \ref{lem:1}. Consequently,\begin{equation}
(q;q)_{\infty}=\sqrt{2\lambda_{n}}\exp\left(\frac{\pi}{24\lambda_{n}}-\frac{\pi\lambda_{n}}{6}\right)\left\{ 1+\mathcal{O}(e^{-4\pi\lambda_{n}})\right\} \label{eq:89}\end{equation}
 \begin{equation}
\frac{1}{(q;q)_{\infty}}=\frac{\exp\left(\frac{\pi\lambda_{n}}{6}-\frac{\pi}{24\lambda_{n}}\right)}{\sqrt{2\lambda_{n}}}\left\{ 1+\mathcal{O}(e^{-4\pi\lambda_{n}})\right\} \label{eq:90}\end{equation}
 as $n\to\infty$. 
\end{proof}
\begin{lem}
\label{lem:3}Given a sequence of positive numbers $\left\{ \lambda_{n}\right\} _{n=1}^{\infty}$,
let \begin{equation}
q=e^{-\pi\lambda_{n}^{-1}},\quad\lim_{n\to\infty}\lambda_{n}=\infty,\quad x>0,\label{eq:91}\end{equation}
 then\begin{equation}
(q^{x};q)_{\infty}=\frac{\sqrt{2}\pi^{1-x}\lambda_{n}^{x-1/2}}{\Gamma(x)\exp(\pi\lambda_{n}/6)}\left\{ 1+\mathcal{O}\left(\lambda_{n}^{-1}\log^{2}\lambda_{n}\right)\right\} \label{eq:92}\end{equation}
 as $n\to\infty$.
\end{lem}
\begin{proof}
From the definition of $\Gamma_{q}(x)$ we have\begin{equation}
(q^{x};q)_{\infty}=\frac{(q;q)_{\infty}(1-q)^{1-x}}{\Gamma_{q}(x)}.\label{eq:93}\end{equation}
Clearly,\begin{equation}
\left(1-q\right)^{1-x}=\left(\frac{\pi}{\lambda_{n}}\right)^{1-x}\left\{ 1+\mathcal{O}\left(\lambda_{n}^{-1}\right)\right\} \label{eq:94}\end{equation}
as $n\to\infty$. In \cite{Zhang3} we proved that for $\Re(z)>0$
\begin{equation}
\Gamma_{q}(z)=\Gamma(z)\left\{ 1+\mathcal{O}\left((1-q)\log^{2}(1-q)\right)\right\} ,\label{eq:95}\end{equation}
 as $q\to1$, thus for $x>0$, \begin{equation}
\Gamma_{q}(x)=\Gamma(x)\left\{ 1+\mathcal{O}\left(\lambda_{n}^{-1}\log^{2}\lambda_{n}\right)\right\} \label{eq:96}\end{equation}
as $n\to\infty$. Combine these equations with equation \eqref{eq:83}
to obtain the equation \eqref{eq:92}.

Take $\lambda=0,\quad\tau=2,\quad m=2n$ in Theorem 2.2 of \cite{Zhang2}
we get the following result:
\end{proof}
\begin{lem}
\label{lem:4}Assume that $z\in\mathbb{C}\backslash\left\{ 0\right\} $,
$\ell>0$ and \eqref{eq:30}, we have \begin{equation}
g(q^{-4n\ell}z;q)=z^{2n}q^{-4n^{2}\ell}\left\{ \theta_{4}\left(z^{-1};q^{\ell}\right)+r_{g}(n|1)\right\} ,\label{eq:97}\end{equation}
 and\begin{align}
|r_{g}(n|1)| & \le\frac{2^{s+r+3}\theta_{3}\left(|z|^{-1};q^{\ell}\right)}{(a_{1},\dotsc,a_{r};q)_{\infty}}\left\{ \frac{q^{n+1}}{1-q}+\frac{q^{\ell n^{2}}}{\left|z\right|^{n}}\right\} \label{eq:98}\end{align}
 for $n$ sufficiently large. In particular we have\begin{align}
{}_{s}\phi_{r}\left(\begin{array}{c}
a_{1},\dotsc,a_{r}\\
b_{1},\dotsc,b_{s}\end{array}\vert q,zq^{-4n\ell}\right) & =\frac{(q,b_{1},\dotsc,b_{s};q)_{\infty}z^{2n}\left\{ \theta_{4}\left(z^{-1}q^{\ell};q^{\ell}\right)+r_{\phi}(n|1)\right\} }{(a_{1},\dotsc,a_{r};q)_{\infty}q^{2\ell n(2n+1)}},\label{eq:99}\end{align}
 and\begin{align}
|r_{\phi}(n|1)| & \le\frac{2^{s+r+3}\theta_{3}\left(|z|^{-1}q^{\ell};q^{\ell}\right)}{(a_{1},\dotsc,a_{r};q)_{\infty}}\left\{ \frac{q^{n+1}}{1-q}+\frac{q^{\ell n^{2}+\ell n}}{\left|z\right|^{n}}\right\} \label{eq:100}\end{align}
 for $n$ sufficiently large, where \begin{equation}
\ell:==\frac{s+1-r}{2}>0.\label{eq:101}\end{equation}

\end{lem}
Similarly, if we take $\lambda=0$ and $\tau=\frac{1}{2}$ in Theorem
2.4 of \cite{Zhang2} we get 

\begin{lem}
\label{lem:5}Assume that $z\in\mathbb{C}\backslash\left\{ 0\right\} $,
$\ell>0$, we have \begin{equation}
h_{n}(zq^{-n\ell};q)=(-z)^{\left\lfloor n/2\right\rfloor }q^{-\ell\left[n^{2}-\chi(n)\right]/4}\left\{ \theta_{4}\left(z^{-1};q^{\ell}\right)+r_{h}(n|1)\right\} ,\label{eq:102}\end{equation}
 and\begin{align}
|r_{h}(n|1)| & \le\frac{2^{s+r+2t+5}\theta_{3}(|z|^{-1};q^{\ell})}{(a_{1},\dotsc,a_{r};q)_{\infty}}\left\{ \frac{q^{\left\lfloor n/4\right\rfloor +1}}{1-q}+|z|^{\left\lfloor n/4\right\rfloor }q^{\ell\left\lfloor n/4\right\rfloor ^{2}}+\frac{q^{\ell\left\lfloor n/4\right\rfloor ^{2}}}{|z|^{\left\lfloor n/4\right\rfloor }}\right\} \label{eq:103}\end{align}
for $n$ sufficiently large.
\end{lem}

\subsection{Proof for Theorem \ref{thm:g-asymptotics-1}}

From \eqref{eq:24} and \eqref{eq:35} to obtain\begin{align}
\theta_{3}(e^{-2\pi v};e^{-\pi\ell\lambda_{n}^{-1}}) & =\theta_{3}\left(vi|\ell\lambda_{n}^{-1}i\right)\label{eq:104}\\
 & =\sqrt{\frac{\lambda_{n}}{\ell}}e^{\pi\ell^{-1}\lambda_{n}v^{2}}\theta_{3}\left(\frac{\lambda_{n}v}{\ell}\mid\frac{\lambda_{n}i}{\ell}\right)\nonumber \\
 & =\sqrt{\frac{\lambda_{n}}{\ell}}e^{\pi\ell^{-1}\lambda_{n}v^{2}}\left\{ 1+\mathcal{O}(e^{-\ell^{-1}\pi\lambda_{n}})\right\} \nonumber \end{align}
as $n\to\infty$, and the $\mathcal{O}$-term is uniform for all $v\in\mathbb{R}$.
Clearly we have\begin{equation}
\frac{q^{n+1}}{1-q}+q^{\ell n^{2}}e^{-2n\pi v}=\mathcal{O}(\lambda_{n}e^{-\pi n\lambda_{n}^{-1}})\label{eq:105}\end{equation}
as $n\to\infty$, and it is uniform for $v$ in any compact subset
of $\mathbb{R}$. From Lemma \ref{lem:3} we have\begin{align}
(q^{\alpha_{1}},\dots,q^{\alpha_{r}};q)_{\infty} & =\frac{2^{r/2}\pi^{r-\sum_{j=1}^{r}\alpha_{j}}\left\{ 1+\mathcal{O}(\lambda_{n}^{-1}\log^{2}\lambda_{n})\right\} }{e^{r\pi\lambda_{n}/6}\lambda_{n}^{r/2-\sum_{j=1}^{r}\alpha_{j}}\prod_{j=1}^{r}\Gamma(\alpha_{j})}\label{eq:106}\end{align}
as $n\to\infty$. Condition \eqref{eq:43} gives\begin{align}
g(-q^{-4n\ell}z;q) & =\sqrt{\frac{\lambda_{n}}{\ell}}\exp\left\{ \frac{\pi\lambda_{n}}{\ell}\left(v+\frac{2n\ell}{\lambda_{n}}\right)^{2}\right\} \left\{ 1+\mathcal{O}(e^{-\ell^{-1}\pi\lambda_{n}})\right\} \label{eq:107}\end{align}
 as $n\to\infty$, and it is uniform for $v$ in any compact subset
of $\mathbb{R}$.

Since \begin{align}
\theta_{4}\left(z^{-1};q^{\ell}\right) & =\theta_{4}\left(e^{-2\pi v};e^{-\ell\pi\lambda_{n}^{-1}}\right)\label{eq:108}\\
 & =\theta_{4}\left(vi\mid\frac{\ell i}{\lambda_{n}}\right)\nonumber \\
 & =\sqrt{\frac{\lambda_{n}}{\ell}}e^{\pi\ell^{-1}\lambda_{n}v^{2}}\theta_{2}\left(\frac{\lambda_{n}v}{\ell}\mid\frac{i\lambda_{n}}{\ell}\right)\nonumber \\
 & =2\sqrt{\frac{\lambda_{n}}{\ell}}\exp\left(\frac{\pi\lambda_{n}v^{2}}{\ell}-\frac{\pi\lambda_{n}}{4\ell}\right)\nonumber \\
 & \times\cos\frac{\pi\lambda_{n}v}{\ell}\left\{ 1+\mathcal{O}\left(e^{-2\pi\ell^{-1}\lambda_{n}}\right)\right\} ,\nonumber \end{align}
 as $n\to\infty$ and its uniformly in $v\in\mathbb{R}$. Thus we
have\begin{align}
g(q^{-4n\ell}z;q) & =\exp\left\{ \frac{\pi\lambda_{n}}{\ell}\left(v+\frac{2n\ell}{\lambda_{n}}\right)^{2}-\frac{\pi\lambda_{n}}{4\ell}\right\} \label{eq:109}\\
 & \times2\sqrt{\frac{\lambda_{n}}{\ell}}\left\{ \cos\frac{\pi\lambda_{n}v}{\ell}+\mathcal{O}\left(e^{-2\ell^{-1}\pi\lambda_{n}}\right)\right\} \nonumber \end{align}
 as $n\to\infty$, and it is uniform on any compact subset of $\mathbb{R}$.

\subsection{Proof for Corollary \ref{cor:ramanujan}}

By Lemma \ref{lem:2} and Theorem \ref{thm:g-asymptotics-1} we have\begin{align}
A_{q}(-q^{-4n}z) & =\frac{g(-;-;q;1;-q^{-4n}z)}{(q;q)_{\infty}}\label{eq:110}\\
 & =\exp\left\{ \pi\lambda_{n}\left(v+\frac{2n}{\lambda_{n}}\right)^{2}+\frac{\pi\lambda_{n}}{6}-\frac{\pi}{24\lambda_{n}}\right\} \nonumber \\
 & \times\frac{1}{\sqrt{2}}\left\{ 1+\mathcal{O}\left(e^{-\pi\lambda_{n}}\right)\right\} ,\nonumber \end{align}
and\begin{eqnarray}
A_{q}(q^{-4n}z) & = & \frac{g(-;-;q;1;q^{-4n}z)}{(q;q)_{\infty}}\label{eq:111}\\
 & = & \exp\left\{ \pi\lambda_{n}\left(v+\frac{2n}{\lambda_{n}}\right)^{2}-\frac{\pi\lambda_{n}}{12}-\frac{\pi}{24\lambda_{n}}\right\} \nonumber \\
 & \times & \sqrt{2}\left\{ \cos\pi\lambda_{n}v+\mathcal{O}\left(e^{-2\pi\lambda_{n}}\right)\right\} \nonumber \end{eqnarray}
 as $n\to\infty$, it is uniform on any compact subset of $\mathbb{R}$.

\subsection{Proof for Corollary \ref{cor:jackson}}

Apply Lemma \ref{lem:2} and Theorem \ref{thm:g-asymptotics-1} to
equation \eqref{eq:14} to get

\begin{align}
J_{\nu}^{(2)}(2i\sqrt{zq^{-\nu}}q^{-2n};q) & =\frac{g(-;q^{\nu+1};q;1;-zq^{-4n})}{(q;q)_{\infty}^{2}(i\sqrt{zq^{-\nu}}q^{-2n})^{-\nu}}\label{eq:112}\\
 & =\frac{\exp\left(\frac{\pi\lambda_{n}}{3}-\frac{\pi}{12\lambda_{n}}+\frac{\nu^{2}\pi}{4\lambda_{n}}+\frac{\nu\pi i}{2}\right)}{2\sqrt{\lambda_{n}}}\nonumber \\
 & \times\exp\left\{ \pi\lambda_{n}\left(v+\frac{4n+\nu}{2\lambda_{n}}\right)^{2}\right\} \left\{ 1+\mathcal{O}(e^{-\pi\lambda_{n}})\right\} ,\nonumber \end{align}
 and \begin{align}
J_{\nu}^{(2)}(2\sqrt{zq^{-\nu}}q^{-2n};q) & =\frac{g(-;q^{\nu+1};q;1;zq^{-4n})}{(q;q)_{\infty}^{2}(\sqrt{zq^{-\nu}}q^{-2n})^{-\nu}}\label{eq:113}\\
 & =\frac{\exp\left(\frac{\pi\lambda_{n}}{12}-\frac{\pi}{12\lambda_{n}}+\frac{\nu^{2}\pi}{4\lambda_{n}}\right)}{\sqrt{\lambda_{n}}}\nonumber \\
 & \times\exp\left\{ \pi\lambda_{n}\left(v+\frac{4n+\nu}{2\lambda_{n}}\right)^{2}\right\} \left\{ \cos\pi\lambda_{n}v+\mathcal{O}\left(e^{-2\pi\lambda_{n}}\right)\right\} \nonumber \end{align}
 as $n\to\infty$, it is uniform on any compact subset of $\mathbb{R}$.

\subsection{Proof for Corollary \ref{cor:confluent}}

Apply Lemma \ref{lem:2} , Lemma \ref{lem:3} and Theorem \ref{thm:g-asymptotics-1}
to equation \eqref{eq:14} to get\begin{align}
{}_{s}\phi_{r}\left(\begin{array}{c}
q^{\alpha_{1}},\dotsc,q^{\alpha_{r}}\\
q^{\beta_{1}},\dotsc,q^{\beta_{s}}\end{array}\vert q,-zq^{-\ell(4n-1)}\right) & =\frac{(q,q^{\beta_{1}},\dots,q^{\beta_{s}};q)_{\infty}g(-q^{-4n\ell}z;q)}{(q^{\alpha_{1}},\dots,q^{\alpha_{r}};q)_{\infty}}\label{eq:114}\\
 & =\frac{(2\pi)^{\ell}\pi^{\rho}\prod_{j=1}^{r}\Gamma(\alpha_{j})}{\sqrt{\ell}\lambda_{n}^{\rho-1/2}\exp\left(\ell\pi\lambda_{n}/3\right)\prod_{j=1}^{s}\Gamma(\beta_{j})}\nonumber \\
 & \times\left\{ \exp\frac{\pi\lambda_{n}}{\ell}\left(v+\frac{2n\ell}{\lambda_{n}}\right)^{2}\right\} \left\{ 1+\mathcal{O}\left(\lambda_{n}^{-1}\log^{2}\lambda_{n}\right)\right\} ,\nonumber \end{align}
 and\begin{align}
{}_{s}\phi_{r}\left(\begin{array}{c}
q^{\alpha_{1}},\dotsc,q^{\alpha_{r}}\\
q^{\beta_{1}},\dotsc,q^{\beta_{s}}\end{array}\vert q,zq^{-\ell(4n-1)}\right) & =\frac{(q,q^{\beta_{1}},\dots,q^{\beta_{s}};q)_{\infty}g(q^{-4n\ell}z;q)}{(q^{\alpha_{1}},\dots,q^{\alpha_{r}};q)_{\infty}}\label{eq:115}\\
 & =\frac{2^{\ell+1}\pi^{\rho+\ell}\prod_{j=1}^{r}\Gamma(\alpha_{j})}{\sqrt{\ell}\lambda_{n}^{\rho-1/2}\exp\left(\frac{\ell\pi\lambda_{n}}{3}+\frac{\pi\lambda_{n}}{4\ell}\right)\prod_{j=1}^{s}\beta_{j}}\nonumber \\
 & \times\left\{ \exp\frac{\pi\lambda_{n}}{\ell}\left(v+\frac{2n\ell}{\lambda_{n}}\right)^{2}\right\} \left\{ \cos\frac{\pi\lambda_{n}v}{\ell}+\mathcal{O}\left(\lambda_{n}^{-1}\log^{2}\lambda_{n}\right)\right\} \nonumber \end{align}

\subsection{Proof for Theorem \ref{thm:h-asymptotics-1}}

Clearly, \begin{align}
\frac{q^{\left\lfloor n/4\right\rfloor +1}}{1-q}+|z|^{\left\lfloor n/4\right\rfloor }q^{\ell\left\lfloor n/4\right\rfloor ^{2}}+\frac{q^{\ell\left\lfloor n/4\right\rfloor ^{2}}}{|z|^{\left\lfloor n/4\right\rfloor }} & =\mathcal{O}\left(e^{-\pi n/(4\lambda_{n})}\right)\label{eq:116}\end{align}
as $n\to\infty$ and it is uniformly on any compact subset of $\mathbb{R}$.
From equations \eqref{eq:42}, \eqref{eq:104}, \eqref{eq:106} and
\eqref{eq:116} we get\begin{align}
h(-zq^{-n\ell};q) & =\exp\left\{ \frac{\pi\lambda_{n}}{\ell}\left(v+\frac{\ell(n-\chi(n))}{2\lambda_{n}}\right)^{2}+\frac{\ell\pi(n-1)\chi(n)}{2\lambda_{n}}\right\} \label{eq:117}\\
 & \times\sqrt{\frac{\lambda_{n}}{\ell}}\left\{ 1+\mathcal{O}(e^{-\ell^{-1}\pi\lambda_{n}})\right\} \nonumber \end{align}
 as $n\to\infty$, it is uniform on any compact subset of $\mathbb{R}$.
Using equations \eqref{eq:42}, \eqref{eq:104}, \eqref{eq:106},
\eqref{eq:108} and \eqref{eq:116} we get\begin{align}
h(zq^{-n\ell};q) & =\exp\left\{ \frac{\pi\lambda_{n}}{\ell}\left(v+\frac{\ell(n-\chi(n))}{2\lambda_{n}}\right)^{2}+\frac{\ell\pi(n-1)\chi(n)}{2\lambda_{n}}-\frac{\pi\lambda_{n}}{4\ell}\right\} \label{eq:118}\\
 & \times2\sqrt{\frac{\lambda_{n}}{\ell}}\left\{ \cos\frac{\pi\lambda_{n}}{\ell}\left(v+\frac{\ell(n-\chi(n))}{2\lambda_{n}}\right)+\mathcal{O}(e^{-2\ell^{-1}\pi\lambda_{n}})\right\} \nonumber \end{align}
 as $n\to\infty$, and it is uniform on any compact subset of $\mathbb{R}$.

\subsection{Proof for Corollary \ref{cor:ismail-masson}}

For $v\in\mathbb{R}$, formula \ref{eq:15}, Lemma \ref{lem:2} and
Theorem \ref{thm:h-asymptotics-1} implies \begin{align}
h_{n}\left(\sinh\pi\left(v+\frac{i}{2}\right)\mid q\right) & =\frac{h(-;-;-;q;1;-e^{2\pi v}q^{-n})}{(-i)^{n}e^{n\pi v}(q;q)_{\infty}}\label{eq:119}\\
 & =\frac{\exp\left\{ \frac{\pi n^{2}}{4\lambda_{n}}+\frac{\pi\lambda_{n}}{6}-\frac{\pi(1+12\chi(n))}{24\lambda_{n}}\right\} }{(-i)^{n}\sqrt{2}}\nonumber \\
 & \times\left\{ \exp\left[\pi\lambda_{n}\left(v-\frac{\chi(n)}{2\lambda_{n}}\right)^{2}\right]\right\} \left\{ 1+\mathcal{O}(e^{-\pi\lambda_{n}})\right\} ,\nonumber \end{align}
 and\begin{align}
h_{n}(\sinh\pi v & \mid q)=\frac{h(-;-;-;q;1;e^{2\pi v}q^{-n})}{(-1)^{n}e^{n\pi v}(q;q)_{\infty}}\label{eq:120}\\
 & =(-1)^{n}\sqrt{2}\exp\left\{ \frac{n^{2}\pi}{4\lambda_{n}}-\frac{(1+12\chi(n))\pi}{24\lambda_{n}}-\frac{\pi\lambda_{n}}{12}\right\} \nonumber \\
 & \times\left\{ \exp\left[\pi\lambda_{n}\left(v-\frac{\chi(n)}{2\lambda_{n}}\right)^{2}\right]\right\} \left\{ \cos\pi\lambda_{n}\left(v+\frac{n-\chi(n)}{2\lambda_{n}}\right)+\mathcal{O}(e^{-2\pi\lambda_{n}})\right\} \nonumber \end{align}
as $n\to\infty$, it is uniform on any compact subset of $\mathbb{R}$.

\subsection{Proof for Corollary \ref{cor:stieltjes-wigert}}

From Lemma\ref{lem:1} and Lemma\ref{lem:2} we have\begin{align}
\frac{1}{(q;q)_{n}(q;q)_{\infty}} & =\frac{1}{(q;q)_{\infty}^{2}}\frac{(q;q)_{\infty}}{(q;q)_{n}}\label{eq:121}\\
 & =\frac{\exp\left\{ \frac{\pi\lambda_{n}}{3}-\frac{\pi}{12\lambda_{n}}\right\} }{2\lambda_{n}}\left\{ 1+\mathcal{O}\left(e^{-4\pi\lambda_{n}}\right)\right\} \nonumber \end{align}
as $n\to\infty$. Hence, equation \eqref{eq:16} and Theorem \ref{thm:h-asymptotics-1}
implies 

\begin{align}
S_{n}(-zq^{-n};q) & =\frac{h(-;-;-;q;1;-zq^{-n})}{(q;q)_{n}(q;q)_{\infty}}\label{eq:122}\\
 & =\frac{\exp\left\{ \frac{\pi\lambda_{n}}{3}+\frac{\pi(n-1)\chi(n)}{2\lambda_{n}}-\frac{\pi}{12\lambda_{n}}\right\} }{2\sqrt{\lambda_{n}}}\nonumber \\
 & \times\left\{ \exp\pi\lambda_{n}\left(v+\frac{n-\chi(n)}{2\lambda_{n}}\right)^{2}\right\} \left\{ 1+\mathcal{O}\left(e^{-\pi\lambda_{n}}\right)\right\} ,\nonumber \end{align}
and\begin{align}
S_{n}(zq^{-n};q) & =\frac{h(-;-;-;q;1;zq^{-n})}{(q;q)_{n}(q;q)_{\infty}}\label{eq:123}\\
 & =\frac{\exp\left\{ \frac{\pi\lambda_{n}}{12}+\frac{\pi(n-1)\chi(n)}{2\lambda_{n}}-\frac{\pi}{12\lambda_{n}}\right\} }{\sqrt{\lambda_{n}}}\nonumber \\
 & \times\left\{ \exp\pi\lambda_{n}\left(v+\frac{n-\chi(n)}{2\lambda_{n}}\right)^{2}\right\} \left\{ \cos\pi\lambda_{n}\left(v+\frac{n-\chi(n)}{2\lambda_{n}}\right)+\mathcal{O}\left(e^{-2\pi\lambda_{n}}\right)\right\} \nonumber \end{align}
as $n\to\infty$, it is uniform on any compact subset of $\mathbb{R}$.

\subsection{Proof for Corollary \ref{cor:laguerre}}

The equations \eqref{eq:17} , \eqref{eq:121} and Theorem \ref{thm:h-asymptotics-1}
gives

\begin{align}
L_{n}^{(\alpha)}(-zq^{-\alpha-n};q) & =\frac{h(-;-;q^{\alpha+1};q;1;-zq^{-n})}{(q;q)_{n}(q;q)_{\infty}}\label{eq:124}\\
 & =\frac{\exp\left\{ \frac{\pi\lambda_{n}}{3}+\frac{\pi(n-1)\chi(n)}{2\lambda_{n}}-\frac{\pi}{12\lambda_{n}}\right\} }{2\sqrt{\lambda_{n}}}\nonumber \\
 & \times\left\{ \exp\pi\lambda_{n}\left(v+\frac{n-\chi(n)}{2\lambda_{n}}\right)^{2}\right\} \left\{ 1+\mathcal{O}\left(e^{-\pi\lambda_{n}}\right)\right\} ,\nonumber \end{align}
 and\begin{align}
L_{n}^{(\alpha)}(zq^{-\alpha-n};q) & =\frac{h(-;-;q^{\alpha+1};q;1;zq^{-n})}{(q;q)_{n}(q;q)_{\infty}}\label{eq:125}\\
 & =\frac{\exp\left\{ \frac{\pi\lambda_{n}}{12}+\frac{\pi(n-1)\chi(n)}{2\lambda_{n}}-\frac{\pi}{12\lambda_{n}}\right\} }{\sqrt{\lambda_{n}}}\nonumber \\
 & \times\left\{ \exp\pi\lambda_{n}\left(v+\frac{n-\chi(n)}{2\lambda_{n}}\right)^{2}\right\} \left\{ \cos\pi\lambda_{n}\left(v+\frac{n-\chi(n)}{2\lambda_{n}}\right)+\mathcal{O}\left(e^{-2\pi\lambda_{n}}\right)\right\} \nonumber \end{align}
as $n\to\infty$, and it is uniform for any compact subset of $\mathbb{R}$.

\end{document}